\documentclass[12pt]{article}

\usepackage{amsfonts,amssymb,amsmath,amsthm,epsfig,euscript}

\newtheorem{definition}{Definition}
\newtheorem{theorem}{Theorem}

\newtheorem{corollary}{Corollary}
\newtheorem{conjecture}{Conjecture}

\long\def\symbolfootnote[#1]#2{\begingroup
\def\thefootnote{\fnsymbol{footnote}}\footnote[#1]{#2}\endgroup}

\newcommand{\sg}{\sigma}

\def\S{\mathcal{S}}
\def\A{\mathcal{A}}
\def\B{\mathcal{B}}

\title{Classifying Descents According to Parity}

\author{
Sergey Kitaev \\
\small Reykjav\'{i}k University\\[-0.8ex]
\small Ofanleiti 2 \\[-0.8ex]
\small IS-103 Reykjav\'{i}k, Iceland\\[-0.8ex]
\small \texttt{sergey@ru.is}
\and
Jeffrey Remmel \\
\small Department of Mathematics\\[-0.8ex]
\small University of California, San Diego\\[-0.8ex]
\small La Jolla, CA 92093-0112. USA\\[-0.8ex]
\small \texttt{remmel@math.ucsd.edu}
}

\date{\small Submitted: Date 1;  Accepted: Date 2;
 Published: Date 3.\\
\small MR Subject Classifications: 05A15}

\begin{document}
\maketitle

\begin{abstract}
\noindent In this paper we refine the well-known permutation
statistic ``descent" by fixing parity of (exactly) one of the
descent's numbers. We provide explicit formulas for the
distribution of these (four) new statistics. We use certain
differential operators to obtain the formulas. Moreover, we
discuss connection of our new statistics to the Genocchi numbers.
We also provide bijective proofs of some of our results. \\

\end{abstract}

\section{Introduction}

The theory of permutation statistics has a long history and has
grown rapidly in the last few decades. The number of {\em
descents} in a permutation $\pi$, denoted by $des(\pi)$, is a
classical statistic. The {\em descent set} of a permutation
$\pi=\pi_1\pi_2\cdots\pi_n$ is the set of indices $i$ for which
$\pi_i>\pi_{i+1}$. This statistics was first studied by
MacMahon~\cite{macmah} almost a hundred years ago, and it still
plays an important role in the field.

{\em Eulerian numbers} $A(n,k)$ count permutations in the
symmetric group $\S_n$ with $k$ descents and they are the
coefficients of the {\em Eulerian polynomials} $A_n(t)$ defined by
$A_n(t)=\sum_{\pi\in\S_n}t^{1+des(\pi)}$. The Eulerian polynomials
satisfy the identity
$$\sum_{k\geq 0}k^nt^k=\frac{A_n(t)}{(1-t)^{n+1}}.$$

For more properties of the Eulerian polynomials see~\cite{comtet}.

By definition, $\pi_i\pi_{i+1}$ is a descent in a permutation
$\pi=\pi_1\pi_2\cdots\pi_n$ if and only if $\pi_i>\pi_{i+1}$. In
this paper we consider refinement of the notion of a descent by
fixing parity of either the first or the second number (but not
both!) of a descent. For example, if we require that the second
number of a descent must be even, then we do not consider at all
descents ending with an odd number. Thus we introduce four new
permutation statistics.

We provide exact formulas for the distribution of our new
statistics. That is, we give two formulas (for even and odd $n$)
for the number of $n$-permutations having exactly $k$ descents of
a chosen type. These formulas for ``beginning with an even number"
(resp. ``ending with an even number," ``beginning with an odd
number," and ``ending with an odd number") case can be found in
Section~\ref{secR} (resp. \ref{secP}, \ref{secQ}, and~\ref{secM}).
In the first three cases above, we use certain differential
operators to derive explicit answers; we then use trivial
bijections of the symmetric group to itself to proceed with the
fourth case.

In Section~\ref{Genocchi} we link our descents to {\em patterns}
in permutations (see~\cite{bona} for an introduction to the
subject). Using certain patterns, which are an alternative
notation for our descents, we provide equivalent definitions of
the {\em Genocchi numbers} both on even and odd permutations. In
Section~\ref{bijections} we state bijective proofs for certain
results on descents according to parity and, in
Section~\ref{conclusions}, we discuss patterns in which
generalized parity considerations are taken into account.

\section{Definitions and notations}

Let $\S_n$ denote the set of permutations of $\{1,2,\ldots, n\}$.
Let $E = \{0, 2, 4, \ldots \}$ and $O = \{1, 3, 5, \ldots \}$
denote the set of even and odd numbers respectively. Given $\sg =
\sg_1\sg_2 \cdots \sg_n \in \S_n$, we define the following, where
$\chi(\sg_1 \in X)$ is $1$ if $\sg_1$ is of type $X$, and it is
$0$ otherwise.

\begin{itemize}
\item $\overleftarrow{Des}_X(\sg) = \{i: \sg_i > \sg_{i+1} \ \& \
\sg_i \in X\}$ and $\overleftarrow{des}_X(\sg) =
|\overleftarrow{Des}_X(\sg)|$ for $X\in\{E,O\}$;

\item $\overrightarrow{Des}_X(\sg) = \{i: \sg_i > \sg_{i+1} \ \& \
\sg_{i+1} \in X\}$ and $\overrightarrow{des}_X(\sg) =
|\overrightarrow{Des}_X(\sg)|$ for $X\in\{E,O\}$;

\item $R_n(x) = \sum_{\sg \in S_n} x^{\overleftarrow{des}_E(\sg)}$
and $P_n(x,z) = \sum_{\sg \in S_n} x^{\overrightarrow{des}_E(\sg)}
z^{\chi(\sg_1 \in E)}$;

\item $M_n(x) = \sum_{\sg \in S_n} x^{\overleftarrow{des}_O(\sg)}$
and $Q_n(x,z) = \sum_{\sg \in S_n} x^{\overrightarrow{des}_O(\sg)}
z^{\chi(\sg_1 \in O)}$;

\item $R_{n}(x) = \sum_{k=0}^{n} R_{k,n} x^k$ and $P_{n}(x,z) =
\sum_{k=0}^{n} \sum_{j=0}^1 P_{j,k,n} z^jx^k$;

\item $M_{n}(x) = \sum_{k=0}^{n} M_{k,n} x^k$ and $Q_{n}(x,z) =
\sum_{k=0}^{n} \sum_{j=0}^1 Q_{j,k,n} z^jx^k$.

\end{itemize}

Thus our goal in this paper is to study the coefficients
$R_{k,n}$, $P_{j,k,n}$, $M_{k,n}$, and $Q_{j,k,n}$ of the
polynomials $R_n(x)$, $P_n(x)$, $M_n(x)$, and $Q_n(x)$
respectively.

Given any permutation $\sg = \sg_1\sg_2 \cdots \sg_n \in \S_n$, we
label the possible positions of where we can insert $n+1$ to get a
permutation in $\S_{n+1}$ from left to right with 0 to $n$, i.e.,
inserting $n+1$ in position 0 means that we insert $n+1$ at the
start of $\sg$ and for $i \geq 1$, inserting $n+1$ in position $i$
means we insert $n+1$ immediately after $\sg_i$. In such a
situation, we let $\sg^{(i)}$ denote the permutation of $\S_{n+1}$
that results by inserting $n+1$ in position $i$.

Let $\sg^c = (n+1 -\sg_1)(n+1 -\sg_2) \cdots (n+1 -\sg_{n})$
denote the {\em complement} of $\sg$. Clearly, if $n$ is odd, then,
for all $i$,  $\sg_i$ and $n+1 -\sg_i$ have the same parity, whereas
they have opposite parity if $n$ is even. The {\em reverse} of
$\sg$ is the permutation $\sg^r=\sg_{n}\sg_{n-1}\cdots\sg_1$.

\section{Beginning with an even number: properties of $R_n(x)$}\label{secR}

Let $\Delta_{2n}$ be the operation which sends $x^k$ to $k x^{k-1} + (2n+1-k)x^k$
and $\Gamma_{2n+1}$ be the operator that sends $x^k$ to $(k+1)x^k +
(2n+1-k) x^{k+1}$. Then we have the following.

\begin{theorem}\label{thm:1}
The polynomials $\{R_n(x)\}_{n \geq 1}$ satisfy the following recursions.

\begin{enumerate}
\item $R_1(x) = 1$ and $R_2(x) = 1+x$,
\item $R_{2n+1}(x) = \Delta_{2n}(R_{2n}(x))$ for $n \geq 1$, and
\item $R_{2n+2}(x) = \Gamma_{2n+1}(R_{2n+1}(x))$ for $n \geq 1$.
\end{enumerate}
\end{theorem}
\begin{proof}
Part 1 is easy to verify by direct computation.

For part 2, suppose $\sg = \sg_1\sg_2 \cdots \sg_{2n} \in \S_{2n}$
and $\overleftarrow{des}_E(\sg) =k$. It is then easy to see that
if we insert $2n+1$ in position $i$ where $i \in
\overleftarrow{Des}_E(\sg)$, then
$\overleftarrow{des}_E(\sg^{(i)}) =k-1$. However, if we insert
$2n+1$ in position $i$ where $i \notin
\overleftarrow{Des}_E(\sg)$, then
$\overleftarrow{des}_E(\sg^{(i)}) =k$.  Thus $\{\sg^{(i)}: i = 0,
\ldots, 2n\}$ gives a contribution of $kx^{k-1} + (2n+1 -k)x^k$ to
$R_{2n+1}(x)$.

For part 3, suppose $\sg = \sg_1\sg_2 \cdots \sg_{2n+1} \in
\S_{2n+1}$ and $\overleftarrow{des}_E(\sg) =k$. It is then easy to
see that if we insert $2n+2$ in position $i$ where $i \in
\overleftarrow{Des}_E(\sg)$ or $i=2n+1$, then
$\overleftarrow{des}_E(\sg^{(i)}) =k$. Similarly if we insert
$2n+2$ in position $i$ where $i \notin \overleftarrow{Des}_E(\sg)
\cup \{2n+1\}$, then $\overleftarrow{des}_E(\sg^{(i)}) =k+1$. Thus
$\{\sg^{(i)}: i = 0, \ldots, 2n+1\}$ gives a contribution of
$(k+1)x^{k} + (2n+1 -k)x^{k+1}$ to $R_{2n+2}(x)$.
\end{proof}

We can express Theorem~\ref{thm:1} in terms of differential
operators:

\begin{corollary} The polynomials $\{R_n(x)\}_{n \geq 1}$
are given by the following
\begin{enumerate}
\item $R_1(x) = 1$, $R_2(x) = 1+x$, and for $n \geq 1$, \item
$R_{2n+1}(x) = (1-x)\frac{d}{dx}R_{2n}(x)+(1+2n)R_{2n}(x)$ and
\item $R_{2n+2}(x) =
x(1-x)\frac{d}{dx}R_{2n+1}(x)+(1+x(1+2n))R_{2n+1}(x)$.
\end{enumerate}
\end{corollary}

This given, we can easily compute some initial values of $R_n(x)$.

\begin{description}

\item $R_1(x) = 1$.
\item $R_2(x) = 1 + x$.
\item $R_3(x) = 4 + 2x$.
\item $R_4(x) = 4 + 16x + 4x^2$.
\item $R_5(x) = 36 + 72x + 12x^2$.
\item $R_6(x) = 36 + 324x + 324x^2 + 36x^3$.
\item $R_7(x) = 576 + 2592x + 1728x^2 + 144x^3$.
\item $R_8(x) = 576 + 9216x + 20736x^2 + 9216x^3 + 576x^4$.
\end{description}

\begin{theorem}\label{thm:2}  We have $R_{0,2n} = R_{n,2n} = (n!)^2$.
\end{theorem}
\begin{proof}
It is easy to see that the theorem holds for $n =1$.

Now suppose that $\sg = \sg_1\sg_2 \cdots \sg_{2n}$ is such that
$\overleftarrow{des}_E(\sg) =0$. Then we can factor any such
permutation into blocks by reading the permutation from left to
right and cutting after each odd number.  For example if $\sg =
1~2~4~5~3~6~7~9~8~10$, then the blocks of $\sg$ would be $1$,
$2~4~5$, $3$, $6~7$,  $9$, and $8~10$.  Since
$\overleftarrow{des}_E(\sg) =0$, there must be a block of even
numbers at the end which must contain the number $2n$  and which
are arranged in increasing order. We call this final block the
$n$-th block. Every other block must end with an odd number $2k+1$
which can be preceded by any subset of even numbers which are less
than $2k+1$ arranged in increasing order. We call such a block the
$k$-th block. It is then easy to see that there are $n!$ ways to
put the even numbers $2, 4, \ldots, 2n$ into the blocks.  That is,
$2n$ must go in the $n$-th block since if we place $2n$ anywhere
but at the end of the permutation, it would contribute to
$\overleftarrow{Des}_E(\sg)$. Then $2(n-1)$ can either go in block
$n-1$ or block $n$. More generally, $2(n-k)$ can go in any blocks
$(n-k), \ldots, n$.  Once we have arranged the even numbers into
blocks, it is easy to see that we can arrange blocks $0, \ldots,
n-1$ in any order and still get a permutation $\sg$ with
$\overleftarrow{des}_E(\sg)= 0$. It thus follows that there are
$(n!)^2$ such permutations. Hence $R_{0,2n} = (n!)^2$.

Now suppose that $\sg = \sg_1\sg_2 \cdots \sg_{2n}$ is such that
$\overleftarrow{des}_E(\sg) =n$. Then, as above, we can factor any
such permutation into blocks by reading the permutation from left
to right and cutting after each odd number.  For example if $\sg =
4~2~1~5~6~3~10~8~7~9$ then the blocks of $\sg$ would be $4~2~1$,
$5$, $6~3$, $10~8~7$, and  $9$.  Since $\overleftarrow{des}_E(\sg)
=n$, each even number must start a descent and hence, unlike the
case where $\overleftarrow{des}_E(\sg) =0$, there can be no even
numbers at the end. Thus the $n$-th block must be empty. It is
also easy to see that if $0 \leq k \leq n-1$ and there are even
numbers in the $k$-th block, i.e. the block that ends with $2k+1$,
then those numbers must all be greater than $2k+1$ and they must
be arranged in decreasing order.

It is then easy to
see that there are $n!$ ways to put the even numbers $2, 4, \ldots, 2n$ into blocks.
That is, $2n$ may go in any of the blocks 0 through $n-1$,
$2(n-1)$ can go in any of the blocks 0 through $n-2$, etc.
After we have partitioned the even numbers into their respective blocks,
we must arrange the even numbers in each block in decreasing order so that there are a total $n!$ ways to
partition the even numbers into the blocks.
Once we have arranged the even numbers into blocks,
it is easy to see that we can arrange blocks $0, \ldots, n-1$ in any order an
still get a permutation $\sg$ with $\overleftarrow{des}_E(\sg)= n$. It thus follows that there are $(n!)^2$ such permutations. Thus $R_{n,2n} = (n!)^2$.
\end{proof}

\begin{theorem}\label{thm:3} We have $R_{k,2n} = \binom{n}{k}^2 (n!)^2$.
\end{theorem}
\begin{proof}
It is easy to see from Theorem~\ref{thm:1} that we have two
following recursions for the coefficients $R_{k,n}$.

\begin{equation}\label{2n:2n+1rec}
R_{k,2n+1} =  (k+1)R_{k+1,2n} +(2n+1-k)R_{k,2n}
\end{equation}

\begin{equation}\label{2n+1:2n+2rec}
R_{k,2n+2} =  (k+1)R_{k,2n+1} +(2n+2-k)R_{k-1,2n+1}
\end{equation}

Iterating the recursions (\ref{2n:2n+1rec}) and (\ref{2n+1:2n+2rec}),
we see that
\begin{eqnarray*}\label{2n:2n+2rec}
R_{k,2n+2} &=& (k+1)((k+1)R_{k+1,2n} +(2n+1-k)R_{k,2n}) \\
&&+ (2n+2-k) ((k)R_{k,2n} +(2n+2-k)R_{k-1,2n})\\
&=& (k+1)^2 R_{k+1,2n} + (2k(2n+1 -k) + 2n +1)R_{k,2n} + (2n+2-k)^2 R_{k-1,2n}
\end{eqnarray*}
Thus if we solve for $R_{k+1,2n}$, we get that
\begin{equation} \label{keyrec}
R_{k+1,2n} = \frac{R_{k,2n+2} - (2k(2n+1 -k) + 2n +1)R_{k,2n} -
(2n+2-k)^2 R_{k-1,2n}}{(k+1)^2}.
\end{equation}

Since we have the initial conditions that $R_{-1,2n} = 0$ and
$R_{0,2n} = (n!)^2$ for all $n \geq 1$, it is a routine
verification using one's favorite computer algebra system that
$R_{k,2n} = \binom{n}{k}^2 (n!)^2$ is the unique solution to
(\ref{keyrec}).
\end{proof}

It follows that we can also get a simple formula for $R_{k,2n+1}$ for
all $n$ and $k$. That is, it immediately follows from (\ref{2n:2n+1rec}) that

\begin{corollary}\label{oddR} One has
$$R_{k,2n+1}=(k+1){\binom{n}{k+1}}^2 (n!)^2 + (2n+1-k) {\binom{n}{k}}^2(n!)^2=\frac{1}{k+1}{\binom{n}{k}}^2((n+1)!)^2.$$
\end{corollary}

\section{Ending with an even number: properties of $P_n(x,z)$}\label{secP}

Let $\Theta_{2n}$ be the operator that sends $z^0x^k$ to $(n+k+1)z^0x^k + (n-k)z^0x^{k+1}$ and sends  $z^1x^k$ to $(n+k+1)z^1x^k + z^0x^{k+1} + (n-k-1)z^1x^{k+1}$.  Let $\Omega_{2n+1}$ be the operator that sends $z^0x^k$ to
$(n+k+1)z^0x^k + z^1x^{k} + (n-k)z^0x^{k+1}$ and sends
$z^1x^k$ to $(n+k+2)z^1x^k + (n-k)z^1x^{k+1}$. Then we have the following.

\begin{theorem}\label{thm:p1} The polynomials $\{P_n(x,z)\}_{n \geq 1}$ satisfy the following recursions.

\begin{enumerate}
\item $P_1(x,z) = 1$ and $P_2(x,z) = 1 + z$.

\item For all $n \geq 1$,
$P_{2n+1}(x,z) = \Theta_{2n}(P_{2n}(x,z))$.

\item For all $n \geq 1$,
$P_{2n+2}(x,z) = \Omega_{2n+1}(P_{2n+1}(x,z))$
\end{enumerate}
\end{theorem}
\begin{proof}
Part 1 is easy to verify by direct computation.

For part 2, suppose $\sg = \sg_1\sg_2 \cdots \sg_{2n} \in
\S_{2n}$, $\overrightarrow{des}_E(\sg) =k$ and $\sg_1 \notin E$.
It is then easy to see that if we insert $2n+1$ in position $i$
where $i \in \overrightarrow{Des}_E(\sg)$, then
$\overrightarrow{des}_E(\sg^{(i)}) =k$. Similarly, if we insert
$2n+1$ in position $i$ where either $i = 2n$ or $\sg_{i+1}$ is
odd, then $\overrightarrow{des}_E(\sg^{(i)}) =k$. However, if we
insert $2n+1$ in position $i$ where $i \notin
\overleftarrow{Des}_E(\sg)$ and $\sg_{i+1}$ is even, then
$\overrightarrow{des}_E(\sg^{(i)}) =k+1$. In every case,
$\sg^{(i)}$ will start with an odd number so that $\{\sg^{(i)}: i
= 0, \ldots, 2n\}$ gives a contribution of $(1+k+n)z^0x^k +
(n-k)z^0x^{k+1}$ to $P_{2n+1}(x,z)$.

Similarly, suppose that  $\sg = \sg_1\sg_2 \cdots \sg_{2n} \in
\S_{2n}$ and $\overrightarrow{des}_E(\sg) =k$ and $\sg_1 \in E$.
It is then easy to see that if we insert $2n+1$ in position $i$
where $i \in \overrightarrow{Des}_E(\sg)$, then
$\overrightarrow{des}_E(\sg^{(i)}) =k$. Similarly, if we insert
$2n+1$ in position $i >0$ where either $i = 2n$ or $\sg_{i+1}$ is
odd, then $\overrightarrow{des}_E(\sg^{(i)}) =k$. If we insert
$2n+1$ in position $i$ where $i > 0$ and $i \notin
\overleftarrow{Des}_E(\sg)$ and $\sg_{i+1}$ is even, then
$\overrightarrow{des}_E(\sg^{(i)}) =k+1$.  In all of these cases,
$\sg^{(i)}$ will start with an even number. However, if we insert
$2n+1$ in position 0, then $\overrightarrow{des}_E(\sg^{(i)})
=k+1$ but $\sg^{(0)}$ will start with an odd number. Thus in this
case, $\{\sg^{(i)}: i = 0, \ldots, 2n\}$ gives a contribution of
$(1+k+n)z^1x^k + (n-k-1)z^1x^{k+1} + z^0x^{k+1}$ to
$P_{2n+1}(x,z)$.

For part 3, suppose $\sg = \sg_1\sg_2 \cdots \sg_{2n+1} \in
\S_{2n+1}$ and $\overrightarrow{des}_E(\sg) =k$ and $\sg_1 \notin
E$. It is then easy to see that if we insert $2n+2$ in position
$i$ where $i \in \overrightarrow{Des}_E(\sg)$, then
$\overrightarrow{des}_E(\sg^{(i)}) =k$. Similarly, if we insert
$2n+2$ in position $i$ where either $i = 2n+1$ or $\sg_{i+1}$ is
odd and $i > 0$, then $\overrightarrow{des}_E(\sg^{(i)}) =k$. If
we insert $2n+2$ in position $i$ where $i \notin
\overleftarrow{Des}_E(\sg)$ and $\sg_{i+1}$ is even, then
$\overrightarrow{des}_E(\sg^{(i)}) =k+1$. In all such  cases,
$\sg^{(i)}$ will start with an odd number. However if we insert
$2n+2$ in position 0, then $\overrightarrow{des}_E(\sg^{(0)}) =k$
but $\sg^{(0)}$ will start with an even number. Thus in this case,
$\{\sg^{(i)}: i = 0, \ldots, 2n+1\}$ gives a contribution of
$(1+k+n)z^0x^k + zx^k+ (n-k)z^0x^{k+1}$ to $P_{2n+2}(x,z)$.

Similarly, suppose that  $\sg = \sg_1\sg_2 \cdots \sg_{2n+1} \in
\S_{2n+1}$ and $\overrightarrow{des}_E(\sg) =k$ and $\sg_1 \in E$.
It is then easy to see that if we insert $2n+2$ in position $i$
where $i \in \overrightarrow{Des}_E(\sg)$, then
$\overrightarrow{des}_E(\sg^{(i)}) =k$. Also, if we insert $2n+2$
in position $i$ where either $i = 2n+1$ or $\sg_{i+1}$ is odd,
then $\overrightarrow{des}_E(\sg^{(i)}) =k$. If we insert $2n+2$
in position $i$ where $i \notin \overleftarrow{Des}_E(\sg)$ and
$\sg_{i+1}$ is even, then $\overrightarrow{des}_E(\sg^{(i)})
=k+1$.  In all of these cases, $\sg^{(i)}$ will start with an even
number. Thus in this case, $\{\sg^{(i)}: i = 0, \ldots, 2n\}$
gives a contribution of $(1+k+n+1)z^1x^k + (n-k)z^1x^{k+1}$ to
$P_{2n+2}(x,z)$.
\end{proof}

We can express Theorem~\ref{thm:p1} in terms of differential
operators:

\begin{corollary} The polynomials $\{P_n(x,z)\}_{n \geq
1}$ are given by the following
\begin{enumerate}
\item $P_1(x,z) = 1$, $P_2(x,z) = 1+z$, and for all $n \geq 1$,
\item $P_{2n+1}(x,z) =
x(1-x)\frac{\partial}{\partial x}P_{2n}(x,z)+x(1-z)\frac{\partial}{\partial z}P_{2n}(x,z)+(1+n(1+x))P_{2n}(x,z)$
and \item $P_{2n+2}(x,z) =
x(1-x)\frac{\partial}{\partial x}P_{2n+1}(x,z)+z(1-z)\frac{\partial}{\partial z}P_{2n+1}(x,z)+(1+z+n(1+x))P_{2n+1}(x,z)$.
\end{enumerate}
\end{corollary}

This given, we can easily compute the first few polynomials $P_n(x,z)$.
\begin{description}
\item $P_1(x,z) = 1$
\item $P_2(x,z) = 1+z$
\item $P_3(x,z) = 2+2z + 2x$
\item $P_4(x,z) = 4+8z+8x+4xz$
\item $P_5(x,z) = 12 +24z +48x + 24xz + 12x^2$
\item $P_6(x,z) = 36 + 108z + 216x + 216xz + 108x^2 + 36x^2z$
\item $P_7(x,z) = 144 + 432z + 1296x + 1296xz + 1296x^2 + 432x^2z + 144x^3$
\item $P_8(x,z) = 576+ 2304z +6912x + 10368xz + 10368x^2 +6192zx^2 +2304x^3 +
576zx^3$
\end{description}

\begin{theorem}\label{thm:p2}  We have
\begin{enumerate}
\item $P_{0,0,2n} = (n!)^2$.

\item $P_{1,0,2n} = (n!) (n+1)! - (n!)^2 = n (n!)^2$.

\item $P_{0,0,2n+1} = (n!) (n+1)!$.

\item $P_{1,0,2n+1} = \left( (n+1)! \right)^2 - (n!) (n+1)!
= n \left( n! (n+1)! \right)$.
\end{enumerate}
\end{theorem}
\begin{proof}
It is easy to see that the theorem holds for $n =1$.

Now suppose that $\sg = \sg_1\sg_2 \cdots \sg_{2n}$ is such that
$\overrightarrow{des}_E(\sg) =0$. Then we can factor any
such permutation into blocks by reading the permutation from right
to left and cutting after each odd number.  For example if $\sg =
1~2~4~5~6~9~7~3~8~10$, then the blocks of $\sg$ would be $1~2~4$,
$5~6$, $3$, $9$,  $7$, and $3~8~10$. Since
$\overrightarrow{des}_E(\sg) =0$, there may be a block of even
numbers at the start which contains only even numbers that are
arranged in increasing order. We call this final block the $0$-th
block. Every other block, when read from left to right, must start  with an odd number $2k-1$ which
can be followed  by any subset of even numbers which are greater
than $2k-1$ arranged in increasing order. We call such a block the
$k$-th block. It is then easy to see that there are $(n+1)!$ ways
to put the even numbers $2, 4, \ldots, 2n$ into the blocks.  That
is, $2n$ can go in any of the blocks $0,1, \ldots, n$; $2(n-1)$
can go in any of the blocks $0, \ldots, n-1$, etc. More generally,
$2(n-k)$ can go in any of blocks $0, \ldots, n-k$.  Once we have
arranged the even numbers into blocks, it is easy to see that we
can arrange blocks $1, \ldots, n$ in any order and still get a
permutation $\sg$ with $\overrightarrow{des}_E(\sg)= 0$. It thus
follows that there are $n!(n+1)!$ such permutations. Thus
$P_{0,0,2n} + P_{1,0,2n} = n!(n+1)!$.

Now if we consider only permutations $\sg = \sg_1\sg_2 \cdots
\sg_{2n}$ such that $\overrightarrow{des}_E(\sg) =0$ and $\sg_1$
is odd, then we cannot put any even numbers in the 0-th block and
hence there are only $n!$ ways to put the even numbers into
blocks. That is, $2n$ can go in any of the blocks $1, \ldots, n$;
$2(n-1)$ can go in any of the blocks $1, \ldots, n-1$, etc. It
follows that $P_{0,0,2n} = (n!)^2$ and that $P_{1,0,2n} = n!(n+1)!
-(n!)^2 = n (n!)^2$.

Now suppose that $\sg = \sg_1\sg_2 \cdots \sg_{2n+1}$ is such that
$\overrightarrow{des}_E(\sg) =0$. Then again we can factor any such
permutation into blocks by reading the permutation from right to
left and cutting after each odd number.  For example if $\sg = 11~
1~2~4~5~6~9~7~3~8~10$, then the blocks of $\sg$ would be $11$,
$1~2~4$, $5~6$, $3$, $9$,  $7$, and $3~8~10$. Since
$\overrightarrow{des}_E(\sg) =0$, there may be a block of even
numbers at the start which contains only even numbers that are
arranged in increasing order. We call this final block the $0$-th
block. Every other block, when read from left to right, must start
  with an odd number $2k-1$ which
can be followed  by any subset of even numbers which are greater
than $2k-1$ arranged in increasing order. We call such a block the
$k$-th block. Notice that this forces the $(n+1)$-st block, i.e.,
the block which ends with $2n+1$ to have no even numbers in it. It
is then easy to see that there are $(n+1)!$ ways to put the even
numbers $2, 4, \ldots, 2n$ into the blocks.  That is, $2n$ can go
in any of the blocks $0,1, \ldots, n$; $2(n-1)$ can go in any of
the blocks $0, \ldots, n-1$, etc. More generally, $2(n-k)$ can go
in any blocks $0, \ldots, n-k$.  Once we have arranged the even
numbers into blocks, it is easy to see that we can arrange blocks
$1, \ldots, n+1$ in any order and still get a permutation $\sg$
with $\overrightarrow{des}_E(\sg)= 0$. It thus follows that there
are $((n+1)!)^2$ such permutations. Thus $P_{0,0,2n+1} +
P_{1,0,2n+1} = ((n+1)!)^2$.

Now if we consider only permutations $\sg = \sg_1\sg_2 \cdots
\sg_{2n+1}$ such that $\overrightarrow{des}_E(\sg) =0$ and $\sg_1$
is odd, then we cannot put any even numbers in the 0-th block and
hence there are only $n!$ ways to put the even numbers into
blocks. That is, $2n$ can go in any of the blocks $1, \ldots, n$;
$2(n-1)$ can go in any of the blocks $1, \ldots, n-1$, etc. It
follows that $P_{0,0,2n} = n!(n+1)!$ and that $P_{1,0,2n} =
((n+1)!)^2 -n!(n+1)! = n(n!(n+1)!)$.
\end{proof}

\begin{theorem}
For all $0 \leq k \leq n$,
\begin{eqnarray*}
P_{0,k,2n+1} &=& P_{0,n-k,2n+1} \  \mbox{and} \\
P_{1,k,2n+1} &=& P_{1,n-k-1,2n+1}
\end{eqnarray*}
\end{theorem}
\begin{proof}

Now suppose that $\sg_1$ is odd. Then $2n+2 -\sg_1$ is odd and,
for any even number $2m$, if $\sg_i =2m$, then $i > 1$ and $i-1
\in \overrightarrow{Des}_E(\sg) \iff i-1 \notin
\overrightarrow{Des}_E(\sg^c)$. It follows that
$\overrightarrow{des}_E(\sg) = k \iff
\overrightarrow{des}_E(\sg^c) =n-k$. Hence the complementation map
shows $P_{0,k,2n+1} = P_{0,n-k,2n+1}$.

Similarly, suppose that $\sg_1$ is even. Then $2n+2 -\sg_1$ is
even. Moreover, there are $n-1$ even numbers $2m$ such that there
exists an $i >1$ such that $\sg_i = 2m$ and, for any such even
number $2m$, $i-1 \in \overrightarrow{Des}_E(\sg) \iff i-1 \notin
\overrightarrow{Des}_E(\sg^c)$. It follows that
$\overrightarrow{des}_E(\sg) = k \iff
\overrightarrow{des}_E(\sg^c) =n-1-k$. Hence the complementation
map shows $P_{1,k,2n+1} = P_{1,n-k-1,2n+1}$.
\end{proof}

Here is a result about the relationships between the polynomials $R_n(x)$ and
the polynomials $P_n(x,z)$
\begin{theorem} \label{R-Prelation}
For all $k$ and $n$,
$R_{k,2n+1} = P_{0,k,2n+1} + P_{1,k,2n+1}$ so that
$R_{2n+1}(x) = P_{2n+1}(x,1)$.
\end{theorem}
\begin{proof}
This result follows by sending $\sg$ to the reverse of $\sg^c$.
That is, if $\sg_i > \sg_{i+1}$ and $\sg_i \in E$, then
$\sg^c_i < \sg^c_{i+1}$ and $\sg^c_i \in E$ and, hence, in
the reverse of $\sg^c$, $\sg^c_i$ will be part of descent that
ends in an even number.
\end{proof}

\begin{theorem}\label{thm:P-values2n} For all $0 \leq k \leq n$,
\begin{equation*}\label{P-values1}
P_{1,k,2n} = \binom{n-1}{k} \binom{n}{k+1} (n!)^2
\end{equation*}
and
\begin{equation*}\label{P-values0}
P_{0,k,2n} = \binom{n-1}{k} \binom{n}{k} (n!)^2
\end{equation*}
\end{theorem}
\begin{proof}
 It follows from Theorem \ref{thm:p1} that we have the following recursions
for the coefficients of $P_{2n+1}(x,z)$ and $P_{2n+2}(x,z)$.
\begin{equation}\label{0prec2n+1}
P_{0,k,2n+1} = (n+k+1)P_{0,k,2n} + P_{1,k-1,2n} +(n-k+1)P_{0,k-1,2n},
\end{equation}

\begin{equation}\label{1prec2n+1}
P_{1,k,2n+1} = (n+k+1)P_{1,k,2n} + (n-k)P_{1,k-1,2n},
\end{equation}

\begin{equation}\label{0prec2n+2}
P_{0,k,2n+2} = (n+k+1)P_{0,k,2n+1} + (n-k+1)P_{0,k-1,2n+1},
\end{equation}
and
\begin{equation*}\label{1prec2n+2}
P_{1,k,2n+2} = (n+k+2)P_{1,k,2n+1} + P_{0,k,2n+1} +(n-k+1)P_{1,k-1,2n+1}.
\end{equation*}

Thus it follows that
\begin{eqnarray}\label{doubrec0}
P_{0,k,2n+2} &=& (n+k+1)\left( (n+k+1)P_{0,k,2n} + P_{1,k-1,2n} +(n-k+1)P_{0,k-1,2n} \right)\nonumber \\
&& + (n-k+1)\left( (n+k)P_{0,k-1,2n} + P_{1,k-2,2n} +(n-k+2)P_{0,k-2,2n} \right)\nonumber\\
&=& (n+k+1)^2 P_{0,k,2n} + (2n^2+3n - 2k^2 +k +1) P_{0,k-1,2n} \nonumber \\
&& + (n-k+1)(n-k+2)P_{0,k-2,2n} \nonumber \\
&& + (n+k+1)P_{1,k-1,2n} + (n-k+1)P_{1,k-2,2n}.
\end{eqnarray}
Similarly,
\begin{eqnarray}\label{doubrec1}
P_{1,k,2n+2} &=& (n+k+2)\left( (n+k+1)P_{1,k,2n} + (n-k)P_{1,k-1,2n} \right) \nonumber \\
&& + (n+k+1)P_{0,k,2n} + P_{1,k-1,2n} +(n-k+1)P_{0,k-1,2n} \nonumber \\
&& + (n-k+1)\left( (n+k)P_{1,k-1,2n} + (n-k+1)P_{1,k-2,2n} \right) \nonumber \\
&=& (n+k+1)(n+k+2) P_{1,k,2n} + (2n^2+3n - 2k^2 - k +1) P_{1,k-1,2n} \nonumber \\
&& + (n-k+1)^2P_{1,k-2,2n} \nonumber \\
&& + (n+k+1)P_{0,k,2n} + (n-k+1)P_{0,k-1,2n}.
\end{eqnarray}

We note that our proof of Theorem \ref{thm:p2} shows that our
formulas hold for all $n$ when $k=0$. It is also easy to check
that our formulas hold when $n=1$ for all $k$.
Next we consider the case when $k=1$.
In this case $P_{0,1-2,2n} = P_{1,1-2,2n} = 0$ for all $n$ by definition so
that the recursion (\ref{doubrec0})  reduces to
\begin{equation}\label{pval01}
P_{0,1,2n+2} = (n+2)^2 P_{0,1,2n} + (2n^2+3n) P_{0,0,2n} + (n+2)P_{1,0,2n}.
\end{equation}
Given that $P_{0,0,2n} = (n!)^2$ and $P_{1,0,2n} = n (n!)^2$ by
Theorem \ref{thm:p2}, assuming by induction that $P_{0,1,2n} =
(n-1)(n)(n!)^2$, and using~(\ref{pval01}), we obtain that
\begin{eqnarray*}
P_{0,1,2n+2} &=& (n+2)^2 (n-1) n (n!)^2 + (2n^2 +3n) (n!)^2 + (n+2) n (n!)^2 \\
&=& n (n+1) ((n+1)!)^2.
\end{eqnarray*}
Thus our formula for $P_{0,1,2n}$ holds by induction.

Similarly the recursion (\ref{doubrec1})  reduces to
\begin{equation}\label{pval011}
P_{1,1,2n+2} = (n+3)(n+2) P_{1,1,2n} + (2n^2+3n-2) P_{1,0,2n} +
(n+2)P_{0,1,2n} + nP_{0,0,2n}.
\end{equation}
Note that $P_{0,0,2n} = (n!)^2$ and $P_{1,0,2n} = n (n!)^2$ by
Theorem \ref{thm:p2} and we just proved that $P_{0,1,2n} =
(n-1)n(n!)^2$. Thus if we assume by induction that $P_{1,1,2n} =
(n-1)\binom{n}{2}(n!)^2$ then using (\ref{pval011}), we obtain
that
\begin{eqnarray*}
P_{1,1,2n+2}&=& (n+3)(n+2) (n-1) \binom{n}{2} (n!)^2 + (2n^2 +3n-2) n (n!)^2 \\
&& + (n+2) (n-1) n (n!)^2 + n(n!)^2 = n \binom{n+1}{2} ((n+1)!)^2.\\
\end{eqnarray*}
Thus our formula for $P_{1,1,2n}$ holds by induction.

We now are in a position to prove the general cases of our
formulas.  That is, we shall prove our formulas hold for all $n$
by induction on $k$. That is, assume our formulas hold for all $n$
and for all $j <k$ and that they also hold for $2n$ and $k$. Then,
using~(\ref{doubrec0}),
\begin{eqnarray}\label{doubrec01}
P_{0,k,2n+2} &=& (n+k+1)^2 \binom{n-1}{k} \binom{n}{k} (n!)^2  +
(2n^2+3n - 2k^2 +k +1) \binom{n-1}{k-1} \binom{n}{k-1} (n!)^2 \nonumber \\
&& + (n-k+1)(n-k+2)\binom{n-1}{k-2} \binom{n}{k-2} (n!)^2 \nonumber \\
&& + (n+k+1)\binom{n-1}{k-1} \binom{n}{k} (n!)^2 +
(n-k+1)\binom{n-1}{k-2} \binom{n}{k-1} (n!)^2.
\end{eqnarray}
Thus we have to show that
\begin{eqnarray}\label{doubrec02}
\binom{n}{k} \binom{n+1}{k} ((n+1)!)^2 &=& \mbox{(the RHS
of~(\ref{doubrec01})).}
\end{eqnarray}
We can clearly divide both sides of (\ref{doubrec02}) by
$(n-1)! n! (n!)^2$ to obtain
\begin{eqnarray*}\label{doubrec03}
\frac{n(n+1)^3}{k!(n-k)!k!(n-k+1)!} &=&
\frac{(n+k+1)^2}{k!(n-k-1)!k!(n-k)!} \nonumber \\
&& + \frac{2n^2+3n-2k^2+k+1}{(k-1)!(n-k)!(k-1)!(n-k+1)!} \nonumber \\
&& + \frac{(n-k+2)(n-k+1)}{(k-2)!(n-k+1)!(k-2)!(n-k+2)!} \nonumber \\
&& + \frac{(n+k+1)}{(k-1)!(n-k)!(k)!(n-k)!} \nonumber \\
&& + \frac{n-k+1}{(k-2)!(n-k+1)!(k-1)!(n-k)!}.
\end{eqnarray*}
We multiply both sides of (\ref{doubrec03}) by $k!(n-k)! k!
(n-k+1)!$ to get a new identity which is easy to check.

Given that we have proved our formula for $P_{0,k,2n}$, we can
prove our formula for $P_{1,k,2n}$ in a similar manner. That is,
we shall prove our formula for $P_{1,k,2n}$ holds for all $n$ by
induction on $k$. Assume our formula hold for all $n$ and for all
$j <k$ and that it also holds for $2n$ and $k$. Then,
using~\ref{doubrec1},

\begin{eqnarray}\label{doubrec11}
P_{1,k,2n+2}
&=& (n+k+1)(n+k+2) \binom{n-1}{k} \binom{n}{k+1} (n!)^2 \nonumber \\
&&  +
(2n^2+3n - 2k^2 - k +1) \binom{n-1}{k-1} \binom{n}{k} (n!)^2 \nonumber \\
&& + (n-k+1)^2 \binom{n-1}{k-2} \binom{n}{k-1} (n!)^2 \nonumber \\
&& + (n+k+1)\binom{n-1}{k} \binom{n}{k} (n!)^2 \nonumber \\
&& + (n-k+1)
\binom{n-1}{k-1} \binom{n}{k-1} (n!)^2
\end{eqnarray}

Then we must show that
\begin{eqnarray}\label{doubrec12}
\binom{n}{k} \binom{n+1}{k+1} ((n+1)!)^2 = \mbox{(the RHS
of~(\ref{doubrec11})).}
\end{eqnarray}

It is easy to see that we can factor out $(n-1)!n!(n!)^2$ from both
side of (\ref{doubrec12}) to obtain that
\begin{eqnarray}\label{doubrec13}
\frac{n(n+1)^3}{k!(n-k)!(k+1)!(n-k)!} &=&
\frac{(n+k+2)(n+k+1)}{k!(n-k-1)!(k+1)!(n-k-1)!} \nonumber \\
&& + \frac{2n^2+3n-2k^2-k+1}{(k-1)!(n-k)!(k)!(n-k)!} \nonumber \\
&& + \frac{(n-k+1)^2}{(k-2)!(n-k+1)!(k-1)!(n-k+1)!} \nonumber \\
&& + \frac{(n+k+1)}{(k)!(n-k-1)!(k-1)!(n-k)!} \nonumber \\
&& + \frac{(n-k+1)}{(k-1)!(n-k)!(k-1)!(n-k+1)!}.
\end{eqnarray}

If we cancel the terms $(n-k+1)$ from both the numerator and
denominator of the third and fifth terms on the RHS of
(\ref{doubrec13}), then it is easy to see that we can multiple
both sides of (\ref{doubrec13}) by $k!(n-k)!(k+1)!(n-k)!$ to
obtain an identity that is easy to check.
\end{proof}

Having found formulas for $P_{j,k,n}$ for the even values of $n$,
we can easily use the recursions (\ref{0prec2n+1}) and
(\ref{1prec2n+1}) to find formulas for $P_{j,k,n}$ for the odd
values of $n$.

\begin{theorem}\label{oddp} For all $0 \leq k \leq n$,
\begin{equation}\label{0pval2n+1}
P_{0,k,2n+1} = (k+1) \binom{n}{k} \binom{n+1}{k+1}(n!)^2 =
(n+1)\binom{n}{k}^2 (n!)^2
\end{equation}
and
\begin{equation*}\label{1pval2n+1}
P_{1,k,2n+1} = \frac{(n+1)(n-k)}{k+1}\binom{n}{k}^2 (n!)^2.
\end{equation*}
\end{theorem}
\begin{proof}
Using the recursion (\ref{0prec2n+1}) and Theorem \ref{thm:P-values2n}, we
have
\begin{eqnarray*}\label{0prec2n+10}
P_{0,k,2n+1} &=& (n+k+1)\binom{n-1}{k}\binom{n}{k}(n!)^2 \nonumber \\
&&+ \binom{n-1}{k-1}\binom{n}{k}(n!)^2 \nonumber \\
&&+ (n-k+1) \binom{n-1}{k-1}\binom{n}{k-1}(n!)^2 \nonumber \\
&=& \binom{n}{k}(n!)^2\left(n \binom{n-1}{k} + (k+1)\binom{n}{k} \right)
\end{eqnarray*}
Here we have used the identities that $(n-k+1)\binom{n}{k-1} = k
\binom{n}{k}$ and $\binom{n-1}{k} + \binom{n-1}{k-1} =
\binom{n}{k}$. It is then easy to verify that
\begin{eqnarray*}
n \binom{n-1}{k} + (k+1)\binom{n}{k} = (n+1) \binom{n}{k}
\end{eqnarray*}
and we get desired.

By Theorem \ref{R-Prelation} and Corollary \ref{oddR}, we
have that

\begin{eqnarray*}\label{RPrelat2n+1}
R_{k,2n+1} &=&  P_{0,k,2n+1} + P_{1,k,2n+1} =
\frac{1}{k+1}\binom{n}{k}^2((n+1)!)^2.
\end{eqnarray*}
Thus, using~(\ref{0pval2n+1}), we get that
\begin{equation*}
P_{1,k,2n+1} = \frac{1}{k+1}\binom{n}{k}^2((n+1)!)^2 -
(n+1)\binom{n}{k}^2 (n!)^2
\end{equation*}
which leads to the result after a simplification.

\end{proof}

As  a corollary to Theorem~\ref{thm:P-values2n}, one has
$P_{0,k,2n} = P_{1,n-1 -k,2n}$, and, using Theorem~\ref{thm:3}, we
get that
$$R_{k,2n} = P_{0,k,2n} + P_{1,k-1,2n}.$$

\section{Ending with an odd number: properties of $Q_n(x,z)$}\label{secQ}

Let $\Phi_{2n}$ be the operator that sends $z^0x^k$ to $z^1x^k +
(n+k)z^0x^k + (n-k)z^0x^{k+1}$ and sends  $z^1x^k$ to
$(n+k+1)z^1x^k + (n-k)z^1x^{k+1}$.  Let $\Psi_{2n+1}$ be the
operator that sends $z^0x^k$ to $(n+k+1)z^0x^k +
(n-k+1)z^0x^{k+1}$ and sends $z^1x^k$ to $(n+k+1)z^1x^k +
z^0x^{k+1} + (n-k)z^1x^{k+1}$. Then we have the following.

\begin{theorem}\label{thm:p11} The polynomials $\{Q_n(x,z)\}_{n \geq 1}$ satisfy the following recursions.
\begin{enumerate}
\item $Q_1(x,z) = z$ and $Q_2(x,z) = z + x$.

\item For all $n \geq 1$, $Q_{2n+1}(x,z) =
\Phi_{2n}(Q_{2n}(x,z))$.

\item For all $n \geq 1$, $Q_{2n+2}(x,z) =
\Psi_{2n+1}(Q_{2n+1}(x,z))$
\end{enumerate}
\end{theorem}
\begin{proof}
Part 1 is easy to verify by direct computation.

For part 2, suppose $\sg = \sg_1\sg_2 \cdots \sg_{2n} \in
\S_{2n}$, $\overrightarrow{des}_O(\sg) =k$ and $\sg_1 \notin O$.
It is then easy to see that if we insert $2n+1$ in position $i$
where $i \in \overrightarrow{Des}_O(\sg)$, then
$\overrightarrow{des}_O(\sg^{(i)}) =k$. Similarly, if we insert
$2n+1$ in position $i$ where either $i = 2n$ or $\sg_{i+1}$ is
even, then $\overrightarrow{des}_O(\sg^{(i)}) =k$. Notice, that
inserting 2n+1 in position 0 produces $\sg^{(0)}$ starting from an
odd number giving a contribution of $zx^k$ rather than $x^k$ to
$Q_{2n+1}(x,z)$. However, if we insert $2n+1$ in position $i$ in
front of remaining odd numbers then
$\overrightarrow{des}_O(\sg^{(i)}) =k+1$. In every case but one,
$\sg^{(i)}$ starts with an even number so that $\{\sg^{(i)}: i =
0, \ldots, 2n\}$ gives a contribution of $z^1x^k +
(k+1+(n-1))z^0x^k + (n-k)z^0x^{k+1}$.

Similarly, suppose that  $\sg = \sg_1\sg_2 \cdots \sg_{2n} \in
\S_{2n}$, $\overrightarrow{des}_O(\sg) =k$ and $\sg_1 \in O$. It
is then easy to see that if we insert $2n+1$ in position $i$ where
$i \in \overrightarrow{Des}_O(\sg)$, then
$\overrightarrow{des}_O(\sg^{(i)}) =k$. Similarly, if we insert
$2n+1$ in position $i$ where either $i = 2n$ or $\sg_{i+1}$ is
even, then $\overrightarrow{des}_O(\sg^{(i)}) =k$. If we insert
$2n+1$ in position $i$ where $i \notin \overleftarrow{Des}_O(\sg)$
and $\sg_{i+1}$ is odd, then $\overrightarrow{des}_O(\sg^{(i)})
=k+1$.  In all of these cases, $\sg^{(i)}$ will start with an odd
number. Thus in this case, $\{\sg^{(i)}: i = 0, \ldots, 2n\}$
gives a contribution of $(n+k+1)z^1x^k + (n-k)z^1x^{k+1}$ to
$Q_{2n+1}(x,z)$.

For part 3, suppose $\sg = \sg_1\sg_2 \cdots \sg_{2n+1} \in
\S_{2n+1}$, $\overrightarrow{des}_O(\sg) =k$ and $\sg_1 \notin O$.
It is then easy to see that if we insert $2n+2$ in position $i$
where $i \in \overrightarrow{Des}_O(\sg)$, then
$\overrightarrow{des}_O(\sg^{(i)}) =k$. Similarly, if we insert
$2n+2$ in position $i$ where either $i = 2n+1$ or $\sg_{i+1}$ is
even, then $\overrightarrow{des}_O(\sg^{(i)}) =k$. If we insert
$2n+2$ in position $i$ where $i \notin \overleftarrow{Des}_O(\sg)$
and $\sg_{i+1}$ is odd, then $\overrightarrow{des}_O(\sg^{(i)})
=k+1$. In all such  cases, $\sg^{(i)}$ will start with an even
number. Thus in this case, $\{\sg^{(i)}: i = 0, \ldots, 2n+1\}$
gives a contribution of $(1+k+n)z^0x^k + (n-k+1)z^0x^{k+1}$ to
$Q_{2n+2}(x,z)$.

Similarly, suppose that  $\sg = \sg_1\sg_2 \cdots \sg_{2n+1} \in
\S_{2n+1}$, $\overrightarrow{des}_O(\sg) =k$ and $\sg_1 \in O$. It
is then easy to see that if we insert $2n+2$ in position $i$ where
$i \in \overrightarrow{Des}_O(\sg)$, then
$\overrightarrow{des}_O(\sg^{(i)}) =k$. Also, if we insert $2n+2$
in position $i$ where either $i = 2n+1$ or $\sg_{i+1}$ is even,
then $\overrightarrow{des}_O(\sg^{(i)}) =k$. If we insert $2n+2$
in position $i>0$ where $i \notin \overleftarrow{Des}_O(\sg)$ and
$\sg_{i+1}$ is odd, then $\overrightarrow{des}_O(\sg^{(i)}) =k+1$.
In all of these cases, $\sg^{(i)}$ will start with an odd number.
However, $\sg^{(0)}$ starts with an even number and
$\overrightarrow{des}_O(\sg^{(0)}) =k+1$. Thus in this case,
$\{\sg^{(i)}: i = 0, \ldots, 2n\}$ gives a contribution of
$(1+k+n)z^1x^k + x^{k+1}+(n-k)z^1x^{k+1}$ to $Q_{2n+2}(x,z)$.
\end{proof}

We can express Theorem~\ref{thm:p11} in terms of differential
operators:

\begin{corollary} The polynomials $\{Q_n(x,z)\}_{n \geq
1}$ are given by the following
\begin{enumerate}
\item $Q_1(x,z) = z$, $Q_2(x,z) = z+x$, and for all $n \geq 1$,
\item $Q_{2n+1}(x,z) =
x(1-x)\frac{\partial}{\partial x}Q_{2n}(x,z)+z(1-z)\frac{\partial}{\partial z}Q_{2n}(x,z)+(z+n(1+x))Q_{2n}(x,z)$\item
$Q_{2n+2}(x,z) =
x(1-x)\frac{\partial}{\partial x}Q_{2n+1}(x,z)+x(1-z)\frac{\partial}{\partial z}Q_{2n+1}(x,z)+(1+n)(1+x)Q_{2n+1}(x,z)$.
\end{enumerate}
\end{corollary}

This given, we can easily compute the first few polynomials
$Q_n(x,z)$.
\begin{description}
\item $Q_1(x,z) = z$ \item $Q_2(x,z) = z+x$ \item $Q_3(x,z) =
2z+2x + 2xz$ \item $Q_4(x,z) = 4z+8x+8zx+4x^2$ \item $Q_5(x,z) =
12z +24x +48xz + 24x^2 + 12x^2z$ \item $Q_6(x,z) = 36z + 108x +
216xz + 216x^2 + 108x^2z + 36x^3$ \item $Q_7(x,z) = 144z + 432x +
1296xz + 1296x^2 + 1296x^2z + 432x^3 + 144x^3z$ \item $Q_8(x,z) =
576z+ 2304x +6912xz + 10368x^2 + 10368x^2z +6192x^3 +2304x^3z +
576x^4$
\end{description}

Forms of the polynomials $P_n(x,z)$ and $Q_n(x,z)$ suggest the
following result.

\begin{theorem}\label{thm:p12} For $n\geq 1$, $Q_n(x,z)=\Xi(P_n(x,z))$, where
$\Xi$ is the operator that sends $z^0x^k$ to $z^1x^k$ and $z^1x^k$
to $x^{k+1}$. In other words, the number of $n$-permutations $\sg$
beginning with an odd number and having
$\overrightarrow{des}_E(\sg) =k$ is equal to that of
$n$-permutations $\pi$ beginning with an odd number and having
$\overrightarrow{des}_O(\pi) =k$; also, the number of
$n$-permutations $\sg$ beginning with an even number and having
$\overrightarrow{des}_E(\sg) =k$ is equal to that of
$n$-permutations $\pi$ beginning with an even number and having
$\overrightarrow{des}_O(\pi) =k+1$.
\end{theorem}

\begin{proof} We prove the result by induction on $n$. The
statement is true for $n=1,2$, which can be seen directly from the
polynomials. Suppose the statement is true for $2n$, that is,
$Q_{2n}(x,z)=\Xi(P_{2n}(x,z))$. Now
$Q_{2n+1}(x,z)=\Phi_{2n}(Q_{2n}(x,z))=\Phi_{2n}(\Xi(P_{2n}(x,z)))$
and we want to show this to be equal to
$\Xi(P_{2n+1}(x,z))=\Xi(\Theta_{2n}(P_{2n}(x,z))).$ In other
words, we want to prove that the operator $\Phi_{2n}(\Xi(\cdot))$
is identical to the operator $\Xi(\Theta_{2n}(\cdot))$ which can
be checked directly by finding the images of $z^0x^k$ and
$z^1x^k$.

Suppose now the statement is true for $2n+1$. One can use
considerations as above involving proving that the operator
$\Psi_{2n+1}(\Xi(\cdot))$ is identical to the operator
$\Xi(\Omega_{2n+1}(\cdot))$, to show that the statement is also
true for $2n+2$.
\end{proof}

The following corollary to Theorem~\ref{thm:p12} is easy to see.
\begin{corollary}\label{Q:values}
For all $0 \leq k \leq n$, $P_{0,k,n}=Q_{1,k,n}$ and
$P_{1,k,n}=Q_{0,k+1,n}$, and thus
$$\begin{array}{l}
Q_{0,k,2n}=\binom{n-1}{k-1} \binom{n}{k} (n!)^2,\\[2mm]
Q_{0,k,2n+1}=\frac{(n+1)(n-k+1)}{k}\binom{n}{k-1}^2 (n!)^2,\\[2mm]
Q_{1,k,2n}=\binom{n-1}{k} \binom{n}{k} (n!)^2,\\[2mm]
Q_{1,k,2n+1}=\binom{n}{k}^2 n!(n+1)!.
\end{array}$$
\end{corollary}

\section{Beginning with an odd number: properties of $M_n(x)$}\label{secM}

\begin{theorem}\label{thm:M-values} For all $0 \leq k \leq n$,
\begin{equation}
M_{k,2n} = \frac{n+1}{k+1}\binom{n-1}{k} \binom{n}{k} (n!)^2 = \binom{n-1}{k}
\binom{n+1}{k+1}(n!)^2
\end{equation}
 and
\begin{equation}
M_{k,2n+1} = \frac{1}{n-k+1}\binom{n}{k}^2((n+1)!)^2 = \binom{n}{k}
\binom{n+1}{k} n! (n+1)!.
\end{equation}
\end{theorem}

\begin{proof}
To prove the first result, note that $M_{k,2n} =
P_{0,k,2n}+P_{1,k,2n}$ since there is a one-to-one correspondence
between permutations counted by $M_{k,2n}$ and the
$2n$-permutations $\sg$ with $\overrightarrow{des}_E(\sg)=k$.
Indeed, a bijection is given by taking the reverse and then the
complement (the parity of the letters will be changed). Now we
simply apply Theorem~\ref{thm:P-values2n} and simplify
$P_{0,k,2n}+P_{1,k,2n}$.

To prove the second result, note that $M_{k,2n+1} =
Q_{0,k,2n+1}+Q_{1,k,2n+1}$ since there is a one-to-one
correspondence between permutations counted by $M_{k,2n+1}$ and
the $(2n+1)$-permutations $\sg$ with
$\overrightarrow{des}_O(\sg)=k$. Indeed, a bijection is given by
taking the reverse and then the complement (the parity of the
letters will be unchanged). Now we simply use
Corollary~\ref{Q:values} and simplify $Q_{0,k,2n+1}+Q_{1,k,2n+1}$.
\end{proof}

\section{Connection to the Genocchi numbers}\label{Genocchi}
Probably the study of {\em Genocchi numbers} goes back to Euler.
The Genocchi numbers can be defined by the following generating
function.
\begin{equation}\label{eq:Genocchi}
\frac{2t}{e^t+1} = t + \sum_{n \geq 1} (-1)^n G_{2n}
\frac{t^{2n}}{(2n)!}.
\end{equation}

These numbers were studied intensively during the last three
decades (see, e.g.,~\cite{ehrste} and references therein). Dumont
\cite{dumont} showed that the
Genocchi number $G_{2n}$ is the number of permutations
$\sg=\sg_1\sg_2 \cdots \sg_{2n+1}$ in $\S_{2n+1}$ such that
$$\begin{array}{ll}
\sg_i<\sg_{i+1} & \mbox{\ if\ } \sg_i \mbox{\ is odd,} \\
\sg_i>\sg_{i+1} & \mbox{\ if\ } \sg_i \mbox{\ is even.}
\end{array}$$
The first few Genocchi numbers are $1, 1, 3, 17, 155, 2073,\ldots$.

Study of distribution of descents is the same as study of
distribution of consecutive occurrences of
the {\em pattern $21$ with no dashes} (see,
e.g.,~\cite{babste} for terminology). Likewise, distribution of
descents according to parity can be viewed as distribution of
consecutive occurrences of certain patterns.

Let us fix some notations. We use $e, o,$ or $*$ as superscripts
for a pattern's letters to require that in an occurrence of the
pattern in a permutation, the corresponding letters must be even,
odd or either. For example, the permutation 25314 has two
occurrences of the pattern $2^*1^o$ (they are 53 and 31, both of
them are occurrences of the pattern $2^o1^o$), one occurrence of
the pattern $1^o2^e$ (namely, 14), no occurrences of the pattern
$1^o2^o$, and no occurrences of the pattern $2^e1^*$.

Given this notation, we can state an alternative definition of the
Genocchi numbers, which follows directly from the definition
above:

\begin{definition}\label{def:1} The Genocchi number $G_{2n}$ is the number of permutations
$\sg=\sg_1\sg_2 \cdots \sg_{2n+1}$ in $\S_{2n+1}$ that avoid
simultaneously the patterns $1^e2^*$ and $2^o1^*$.\end{definition}

Our terminology allows to define the Genocchi numbers on even
permutations as well, due to the following result which we state
as a conjecture (note that in the conjecture only descents
according to parity are involved unlike Definition~\ref{def:1}):

\begin{conjecture}\label{thm:genocchi1} The number of permutations in $\S_{2n}$ that avoid simultaneously
the patterns $2^*1^e$ and $2^e1^*$ is given by the Genocchi number
$G_{2n}$. \end{conjecture}

The following theorem provides two more alternative definitions of
the Genocchi numbers both on even and on odd permutations.
However, the first of these definitions relies on truth of
Conjecture~\ref{thm:genocchi1}.

\begin{theorem} For $n\geq 1$, the number of permutations in $\S_{2n-1}$ that avoid simultaneously
the patterns $2^*1^e$ and $2^e1^*$ is given by $G_{2n}$. Also, the
number of permutations in $\S_{2n}$ that avoid simultaneously the
patterns $1^e2^*$ and $2^o1^*$ is given by $2G_{2n}$.
\end{theorem}

\begin{proof} The first part of the statement follows from Conjecture~\ref{thm:genocchi1}
using an observation that the number $2n$ must be at the end of a
permutation avoiding $2^*1^e$ and $2^e1^*$, and adjoining this
number from the right to any ``good" $(2n-1)$-permutation gives a
``good" $(2n)$-permutation.

For the second part, we note the following: the number $2n+1$ must
be the rightmost number of a $(2n+1)$-permutation avoiding
$1^e2^*$ and $2^o1^*$. Clearly, if we remove this number from a
``good" permutation, we get a ``good" $(2n)$-permutation. The
reverse is not true, since if a ``good" $(2n)$-permutation ends
with an even number, then adjoining $2n+1$ to the right of it
gives an occurrence of $1^e2^*$, whereas adjoining $2n+1$ to the
right of an odd number does not lead to an occurrence of a
prohibited pattern. Thus to prove the statement, we need to prove
that among $1^e2^*-$ and $2^o1^*-$avoiding $(2n)$-permutations
half end with an even number, and half end with an odd number.
This, however, is easy to see since a bijection between these
objects is given by the complement. Indeed, the complement changes
parity (in particular parity of the rightmost number), and also
$\sg$ avoids $1^e2^*$ and $2^o1^*$ if and only if its complement
$\sg^c$, avoids $1^e2^*$ and $2^o1^*$.
\end{proof}

To summarize the section we say that certain descents according to
parity do not only provide alternative definitions for the
Genocchi numbers but also generalize them in sense that instead of
considering (multi-)avoidance of the descents, which gives the
Genocchi numbers, one may consider, for example, their (join)
distribution.

\section{Bijective proofs related to the context}\label{bijections}
It is always nice to be able to provide a bijective proof for an
identity. From the form of coefficients of some polynomials in
this paper, one can see several relations between different groups
of permutations. In this section we provide bijective solutions
for five such relations.

For Subsections~\ref{bij02} and~\ref{bij03}, recall that if $\sg$
is a permutation of $\S_n$ then $\sg^{(i)}$ denotes the
permutation of $\S_{n+1}$ that results by inserting $n+1$ in
position $i$.

\subsection{Bijective proof for the symmetry of
$R_{2n}(x)$.}\label{bijR(x)} Given a permutation $\sg = \sg_1\sg_2
\cdots \sg_{2n} \in \S_{2n}$ with $\overleftarrow{des}_E(\sg) =k$,
apply the complement to the permutation $\sg'=\sg (2n+1)$, that
is, $\sg'$ is obtained from $\sg$ by adding a dummy number
$(2n+1)$ at the end. In the obtained permutation $\sg'^c$, make a
cyclic shift to the left to make the number $(2n+1)$ be the first
one. Remove $(2n+1)$ to get a $2n$-permutation $\sg^*$ with
$\overleftarrow{des}_E(\sg^*) =n-k$. To reverse this procedure,
adjoin $(2n+1)$ from the left to a given permutation $\sg^*$ with
$\overleftarrow{des}_E(\sg^*) =n-k$. Then make a cyclic shift to
the right to make $1$ be the rightmost number. Use the complement
and remove $(2n+1)$ from the obtained permutation to get a
permutation $\sg$ with $\overleftarrow{des}_E(\sg) =k$.

The map described above and its reverse are clearly injective. We
only need to justify that given $k$ occurrences of the descents in
$\sg$, we get $(n-k)$ occurrences in $\sg^*$ (the reverse to this
statement will follow using the same arguments). Notice that
adding $(2n+1)$ at the end does not increase the number of
descents. Since $\sg'$ ends with an odd number, $\sg'^c$ has $n-k$
descents. If $\sg'=A1B(2n+1)$, where $A$ and $B$ are some factors,
then $\sg'^c=A^c(2n+1)B^c1$ and $\sg^*$ is $(2n+1)B^c1A^c$ without
$(2n+1)$. The last thing to observe is that moving $A^c$ to the
end of $\sg'^c$ does not create a new descent since it cannot
start with $1$, also we do not lose any descents since none of
them can end with $(2n+1)$. So, $\sg^*$ has $n-k$ descents.

\subsection{Bijective proof for $R_{k,2n} = P_{0,k,2n} + P_{1,k-1,2n}$. }

A similar solution as that in Subsection~\ref{bijR(x)} works, with
the main difference that we now keep track of whether this is an
odd or even number to the left of 1 in $\sg= \sg_1\sg_2 \cdots
\sg_{2n} \in \S_{2n}$ with $\overleftarrow{des}_E(\sg) =k$. We do
not provide all the justifications in our explanation of the map,
since they are similar to that in Subsection~\ref{bijR(x)}.

Suppose $\sg'=\sg (2n+1)=Ax1B(2n+1)$, where $A$ and $B$ are some
factors and $x$ is a number. Apply the complement {\em and
reverse} to $\sg'$ to get $(\sg')^{cr}=1B^{cr}(2n+1)x^cA^{cr}$.
Make a cyclic shift to the left in $(\sg')^{cr}$ to make the
number $(2n+1)$ be the first one and to get
$\sg^*=(2n+1)x^cA^{cr}1B^{cr}$. One can check that
$\overleftarrow{des}_E(\sg) = \overrightarrow{des}_E(\sg^*) = k$.
Also, parity of $x$ is the same as parity of $x^c$. Now, if we
remove $(2n+1)$ from $\sg^*$ and $x^c$ is even, we loose one
descent obtaining a permutation counted by $P_{1,k-1,2n}$; if we
remove $(2n+1)$ from $\sg^*$ and $x^c$ is odd, the number of
descents in the obtained permutation is the same, $k$, and thus we
get a permutation counted by $P_{0,k,2n}$. Note that if $Ax$ is
the empty word, that is, $\sg$ starts with 1, then this case is
treated as the case ``$x$ is odd" since $\sg^*$ will start with
$(2n+1)1$. Thus one may think of $2n+1$ as the (cyclic)
predecessor of 1 in this case, that is, $x=2n+1$.

The reverse to the map described is easy to see.

\subsection{Bijective proof for $P_{0,k,2n}=P_{1,n-1-k,2n}$.}

This identity states that the number of $2n$-permutations $\sg$
beginning with an odd number and having
$\overrightarrow{des}_E(\sg) =k$ is equal to that of
$2n$-permutations $\pi$ beginning with an even number and having
$\overrightarrow{des}_E(\pi) =n-1-k$. Let $\A_{2n}(k)$ (resp.
$\A_{2n}(n-1-k)$) denote the set of permutations of the first
(resp. second) kind.

Suppose $\pi=x\pi'\in\A_{2n}(n-1-k)$ where $x$ is an even number.
Adjoin a dummy number $2n+1$ to $\pi$ from the right to get the
permutation $\pi_1=(2n+1)x\pi'$ with
$\overrightarrow{des}_E(\pi_1) =n-k$. Apply the compliment to get
the $(2n+1)$-permutation $\pi_2=\pi_1^c=1x^c(\pi')^c$. Since
parity of the numbers preserved after applying the complement,
clearly $\overrightarrow{des}_E(\pi_2) =k$ and $x^c$ is an even
number. Now remove 1 from $\pi_2$ and decrease each number of
$\pi_2$ by 1 to get the $2n$-permutation $\pi_3=y\pi''$ where
$y=x^c-1$ is an odd number and $\overrightarrow{des}_O(\pi_3) =k$
(each descent ending with an even number becomes a descent ending
with an odd number). Using the notation from
Subsection~\ref{bij02}, $\pi_3\in\B_{2n}(k)$, which is the set of
$2n$-permutations $\tau$ beginning with an odd number and having
$\overrightarrow{des}_O(\tau) =k$.

Since all the steps made above are invertible, we have a bijection
between $\A_{2n}(n-1-k)$ and $\B_{2n}(k)$. It now remains to apply
the (recursive) bijection between $\A_{2n}(k)$ and $\B_{2n}(k)$
provided in Subsection~\ref{bij02}, to get the desired bijection
between $\A_{2n}(k)$ and $\A_{2n}(n-1-k)$.

\subsection{Bijective (recursive) proof for
$P_{0,k,n}=Q_{1,k,n}$.}\label{bij02}

This identity states that the number of $n$-permutations $\sg$
beginning with an odd number and having
$\overrightarrow{des}_E(\sg) =k$ is equal to that of
$n$-permutations $\pi$ beginning with an odd number and having
$\overrightarrow{des}_O(\pi) =k$. Let $\A_n(k)$ (resp. $\B_n(k)$)
denote the set of permutations of the first (resp. second) kind.
Of course, $\cup_{k\geq 0}\A_n(k)=\cup_{k\geq 0}\B_n(k)$ is the
set of all permutations in $\S_n$ that begin with an odd number.

For $\S_1$ we map $1\in \A_1(0)$ to $1\in \B_1(0)$, and for $\S_2$
we map $12\in \A_2(0)$ to $12\in \B_2(0)$ (the permutation
$21\in\S_2$ does not start with an odd number, so we do not
consider it).

Now suppose for $\S_{2n}$ we have a bijective map between
$\A_{2n}(k)$ and $\B_{2n}(k)$ for $k=0,1,\ldots,n-1$, and suppose
that a permutation $\sg\in\A_{2n}(k)$ corresponds to a permutation
$\pi\in\B_{2n}(k)$. Based on $\sg$ and $\pi$ we will match $n+k+1$
permutations from $\A_{2n+1}(k)$ to $n+k+1$ permutations from
$\B_{2n+1}(k)$ (case A1 below), as well as $n-k$ permutations from
$\A_{2n+1}(k+1)$ to $n-k$ permutations from $\B_{2n+1}(k+1)$ (case
A2 below). One can see that all the maps below are bijective and
all the permutations we deal with start with an odd number.
\begin{itemize}
\item[A1:]

\begin{itemize} \item[1.1:] Map $\sg^{(2n)}$ to
$\pi^{(2n)}$. \item[1.2:] Map $\sg^{(i)}$ to $\pi^{(j)}$ if
$i$ is the $m$-th descent in $\sg$ and $j$ is the $m$-th descent
in $\pi$ ($m=1,2,\ldots,k$). \item[1.3:] Map $\sg^{(i)}$ to
$\pi^{(j)}$ if $\sg_{i+1}$ is the $m$-th odd number in $\sg$ and
$\pi_{j+1}$ is the $m$-th even number in $\pi$ ($m=1,2,\ldots,n$).
\end{itemize}

\item[A2:] Map $\sg^{(i)}$ to $\pi^{(j)}$ if $\sg_{i+1}$ is the
$m$-th even number such that $i\not\in\overrightarrow{Des}_E(\sg)$
and $\pi_{j+1}$ is the $m$-th odd number such that
$j\not\in\overrightarrow{Des}_O(\pi)$ ($m=1,2,\ldots,n-k$).
\end{itemize}

Note that cases A1 and A2 cover all possible insertions of
$2n+1$ in $\sg$ and $\pi$. \\

Now suppose for $\S_{2n+1}$ we have a bijective map between
$\A_{2n+1}(k)$ and $\B_{2n+1}(k)$ for $k=0,1,\ldots,n$, and
suppose that a permutation $\sg\in\A_{2n+1}(k)$ corresponds to a
permutation $\pi\in\B_{2n+1}(k)$. Based on $\sg$ and $\pi$ we will
match $n+k+1$ permutations from $\A_{2n+1}(k)$ to $n+k+1$
permutations from $\B_{2n+1}(k)$ (case B1 below), as well as $n-k$
permutations from $\A_{2n+1}(k+1)$ to $n-k$ permutations from
$\B_{2n+1}(k+1)$ (case B2 below).

\begin{itemize}
\item[B1:]

\begin{itemize} \item[1.1:] Map $\sg^{(2n+1)}$ to
$\pi^{(2n+1)}$. \item[1.2:] Map $\sg^{(i)}$ to $\pi^{(j)}$ if
$i$ is the $m$-th descent in $\sg$ and $j$ is the $m$-th descent
in $\pi$ ($m=1,2,\ldots,k$). \item[1.3:] Map $\sg^{(i)}$ to
$\pi^{(j)}$ if $\sg_{i+1}$ is the $(m+1)$-st odd number in $\sg$
and $\pi_{j+1}$ is the $m$-th even number in $\pi$
($m=1,2,\ldots,n$).
\end{itemize}

\item[B2:] Map $\sg^{(i)}$ to $\pi^{(j)}$ if $\sg_{i+1}$ is the
$m$-th even number such that $i\not\in\overrightarrow{Des}_E(\sg)$
and $\pi_{j+1}$ is the $(m+1)$-st odd number such that
$j\not\in\overrightarrow{Des}_O(\pi)$ ($m=1,2,\ldots,n-k$).
\end{itemize}

Note that cases B1 and B2 cover all possible insertions of $2n+2$
in $\sg$ and $\pi$, and this finishes the construction of our
bijective map.

\subsection{Bijective (recursive) proof for
$P_{1,k,n}=Q_{0,k+1,n}$.}\label{bij03}

The bijection we provide in this subsection is based on the
(recursive) bijection (let us denote it by $\alpha$) for
$P_{0,k,n}=Q_{1,k,n}$ from Subsection~\ref{bij02}. Using
$n$-permutations $\sg$ and $\pi$, where $\sg=\alpha (\pi)$, $\sg$
and $\pi$ start with odd numbers, and $\overrightarrow{des}_E(\sg)
=\overrightarrow{des}_O(\pi)=k$, we will match $k$ permutations
counted by $P_{1,k-1,n+1}$ with $k$ permutations counted by
$Q_{0,k,n+1}$, as well as $n-k$ permutations counted by
$P_{1,k,n+1}$ with $n-k$ permutations counted by $Q_{0,k+1,n+1}$.
It will be clear that following our procedure, different $\sg$ and
$\pi$ such that $\sg=\alpha (\pi)$, produce different pairs of
permutations matched. Also, it will be easy to see that each
$(n+1)$-permutation starting with an even number will be taken
into account.

Given $\sg$ and $\pi$ with the properties as above, we consider
$2n$ permutations of $\{0,1,\ldots,n\}$ obtained by placing the
(even) number 0 in $\sg$ and $\pi$ in different positions but
position 0. Inserting 0 in $\sg$ in position $i$, where
$i\in\overrightarrow{Des}_E(\sg)$, produces $\sg'$ with
$\overrightarrow{des}_E(\sg')=k$. Inserting 0 in $\pi$ in position
$i$, where $i\in\overrightarrow{Des}_O(\pi)$, produces $\pi'$ with
$\overrightarrow{des}_O(\pi')=k-1$. We match such $\sg'$ and
$\pi'$ if the insertion was inside the $m$-th descent in both
cases ($m=1,2,\ldots,k$).

Inserting 0 in a non-descent position in $\sg$ (resp. $\pi$) gives
$\sg'$ (resp. $\pi'$) with $\overrightarrow{des}_E(\sg')=k+1$
(resp. $\overrightarrow{des}_O(\pi')=k$). We match such $\sg'$ and
$\pi'$ if the insertion in them was in the same (non-descent)
position $m$ counting from left to right ($m=1,2,\ldots,n-k$).

We now increase each number in all the $\sg'$ and $\pi'$
considered above by 1 (evens become odds and vice versa) to get
permutations counted by $Q_{0,k,n+1}$ been matched to permutations
counted by $P_{1,k-1,n+1}$, as well as permutations counted by
$Q_{0,k+1,n+1}$ been matched to permutations counted by
$P_{1,k,n+1}$.

To summarize, our bijection works as follows. Given an
$(n+1)$-permutation $\delta$ starting with an even number and
having, say, $\overrightarrow{des}_E(\delta) =k$, we decrease each
number of $\delta$ by 1 to get $\delta'$, a permutation of
$\{0,1,\ldots,n\}$ which starts with an odd number. If we now
ignore 0 in $\delta'$, then we get an $n$-permutation $\delta''$
starting with an odd number. We then let $\tau''=\alpha
(\delta'')$. Depending on the position of 0 in $\delta'$ we insert
0 into $\tau''$ to get a permutation $\tau'$ of $\{0,1,\ldots,n\}$
which starts with an odd number and has the right properties (see
instructions above). Finally, we increase all the numbers in
$\tau'$ by 1 to get $\tau$ starting with an even number and having
$\overrightarrow{des}_O(\tau) =k+1$. Thus, $\tau$ corresponds to
$\delta$. The reverse to this is similar.

\section{Conclusions}\label{conclusions}
This paper can be viewed as a first step in a more general program
which is to study pattern avoiding conditions permutations where
generalized parity considerations are taking into account. Below
we provide a possible parity generalization.

For any sequence of distinct numbers, $i_1 \cdots i_m$, we let
$red(i_1 \cdots i_m)$ denote the permutation of $\S_m$ whose
elements have the same relative order as $i_1 \cdots i_m$. For
example, $red(5~2~7~8) = 2~1~3~4$. Then given a permutation $\tau
= \tau_1 \cdots \tau_m \in \S_m$, we say that the permutation
$\sigma = \sigma_1 \cdots \sigma_n \in \S_n$ is $\tau$-avoiding if
there is no subsequence $\sigma_{i_1} \sigma_{i_2} \cdots
\sigma_{i_m}$ of $\sigma$ such that $red(\sigma_{i_1} \sigma_{i_2}
\cdots \sigma_{i_m}) = \tau$. Similarly, we say that $\sigma$ has
a $\tau$-match if there is a consecutive subsequence $\sigma_{i_1}
\sigma_{i_1+1} \cdots \sigma_{i_1+m-1}$ such that
$red(\sigma_{i_1} \sigma_{i_1+1} \cdots \sigma_{i_1+m-1}) = \tau$.
There have been many papers in the literature that have studied
the number of $\tau$-avoiding permutations of $\S_n$ or the
distribution of $\tau$-matches for $\S_n$ (e.g., see
\cite{bona,ElizNoy,KM} and references therein). Now we can
generalize the notion of $\tau$-avoiding permutations or
$\tau$-matches by adding parity type conditions. For example, for
any integer $k \geq 2$, we say that a permutation is
parity-$k$-$\tau$-avoiding if there is no subsequence
$\sigma_{i_1} \sigma_{i_2} \cdots \sigma_{i_m}$ of $\sigma$ such
that $red(\sigma_{i_1} \sigma_{i_2} \cdots \sigma_{i_m}) = \tau$
and for all $j$, $\sigma_{i_j} \equiv \tau_j \mod\ k$. Similarly,
we say that $\sigma$ has a parity-$k$-$\tau$-match if there is a
consecutive subsequence $\sigma_{i_1} \sigma_{i_1+1} \cdots
\sigma_{i_1+m-1}$ such that $red(\sigma_{i_1} \sigma_{i_1+1}
\cdots \sigma_{i_1+m-1})=\tau$ and for all $j$, $\sigma_{i_j}
\equiv \tau_j\mod k$. For example, the permutation $\sigma =
3~2~4~5~1$ is not $2~1$-avoiding and has two $2~1$-matches.
However it is parity-$2$-$2~1$ avoiding and, hence, it has no
parity-$2$-$2~1$-matches. Similarly, $\sigma$ has parity-$3$-$2~1$
match since $red(5~1) = 2~1$, $5 \equiv 2 \mod 3$, and $1 \equiv 1
\mod 3$.

We have started to study such generalized parity matching type
conditions. For example, in \cite{KR}, we have generalized the
results of this paper to classify descents according to
equivalence mod $k$ for $k \geq 3$. In \cite{LR}, Liese and Remmel
have studied the distribution of parity-$k$-$\tau$ matches for
$\tau \in \S_2$.

\end{document}